\documentclass{gtart}
\gtart

\usepackage{amsthm}
\usepackage{amsgen}
\usepackage{amsmath,amssymb}
\usepackage{curvesls}
\usepackage{epic}
\usepackage[all]{xy}
\usepackage[dvips]{graphicx}
\usepackage{psfrag}
\usepackage[T1]{fontenc}
\usepackage[sc]{mathpazo}
\newtheorem{thm}{Theorem}[section]
\newtheorem{lem}[thm]{Lemma}
\newtheorem{cor}[thm]{Corollary}
\newtheorem{prop}[thm]{Proposition}

\newtheorem{prob}[thm]{Problem}

\theoremstyle{definition}
\newtheorem{defn}[thm]{Definition}

\newcommand{\Real}{{\mathbb R}}

\newcommand{\Zed}{{\mathbb Z}}

\begin{document}

\title{Triangulating a Cappell-Shaneson knot complement}

\authors{Ryan Budney \\ Benjamin A. Burton \\ Jonathan Hillman}
\addresses{Mathematics and Statistics, University of Victoria, Canada \\
School of Mathematics and Physics,  University of Queensland, Australia \\
School of Mathematics and Statistics, University of Sydney, Australia}
\emails{rybu@uvic.ca \\ bab@maths.uq.edu.au \\ jonathan.hillman@sydney.edu.au}

\begin{abstract} 
We show that one of the Cappell-Shaneson knot complements admits an extraordinarily small 
topological triangulation, containing only two $4$-dimensional simplices.
\end{abstract}


\maketitle

\section{Introduction}\label{INTRODUCTION}

Thurston observed \cite{Thu1} that the figure-8 knot complement in $S^3$ admits a particularly simple triangulation, having two tetrahedra.  If $K \subset S^3$ is the figure-8 knot, the `triangulation' of the complement is what is known as an {\it ideal triangulation,} meaning it is really a triangulation of the $3$-sphere with the knot $K$ crushed to a point or equivalently this is the one-point compactification of the knot complement $S^3 \setminus K$.  The purpose of this note is to describe a similarly simple ideal triangulation of the complement of an embedded $S^2$ in $S^4$.  It turns out this $2$-knot is a Cappell-Shaneson knot \cite{CS}.  Section \ref{primarysec} describes the triangulation and why it is homeomorphic to the complement of a Cappell-Shaneson knot.  

We use the word `triangulation' in a slightly more flexible form than `geometric realization of a simplicial complex' or even `delta complex' in this paper.  In short, a {\it triangulation} in this paper is a space constructed by gluing simplices together via affine identifications of their boundary facets, and where we demand that the characteristic map of each simplex is an embedding when restricted to the interior of the simplex.  

Precisely, the triangulations we use are called `unordered delta complexes' or `generalized triangulations'.  Denote the $n$-simplex by $\Delta^n = \{(x_0,\cdots,x_n) \in \Real^{n+1} : x_i \geq 0 \ \forall i \text{ and } x_0+x_1+\cdots+x_n = 1\}$.  Given $i \in \{0,1,\cdots,n\}$ the $i$-th face map of $\Delta^n$ is the map $f_i : \Delta^{n-1} \to \Delta^n$ given by $f_i(x_0,\cdots,x_{n-1}) = (x_0,x_1,\cdots,x_{i-1}, 0, x_i, x_{i+1}, \cdots, x_{n-1})$. Given a permutation $\sigma \in \Sigma(\{0,1,\cdots,n\})$, the induced automorphism of $\Delta^n$ is given by $\sigma_* : \Delta^n \to \Delta^n$, $\sigma_*(x_0,x_1,\cdots,x_n) = (x_{\sigma(1)}, x_{\sigma(2)}, \cdots, x_{\sigma(n)})$.  An {\it unordered delta complex} is a CW-complex $X$ such that the domains of the attaching maps are the boundaries of simplices (rather than discs), $\phi : \partial \Delta^n \to X^{(n-1)}$, and for each $i$, the composite $\phi \circ f_i = \Phi \circ \sigma_*$ where $\Phi : \Delta^{n-1} \to X^{(n-1)}$ is a characteristic map of the $(n-1)$-skeleton, and $\sigma \in \Sigma(\{0,1,\cdots,n-1\})$ is some permutation. If $\sigma$ is always the identity permutation, this would be an {\it ordered delta complex}. 

We will describe an ideal triangulation of the complement of a knotted $S^2$ in $S^4$.  Precisely, this is an unordered delta complex structure on the one-point compactification of the knot complement, where the only $0$-cell is the point at infinity. In $3$-manifold topology, ideal triangulations turn out to be a useful objects for describing hyperbolic structures on $3$-manifolds \cite{Weeks}. 

The first two authors would like to acknowledge Toshitake Kohno and IPMU (U.Tokyo) for their hospitality.  The authors would like to thank Allen Hatcher and Stephan Tillman for their comments on the first draft of the paper.  The second author is supported by the Australian Research Council under the Discovery Projects funding scheme (project DP1094516), the first author by an NSERC Discovery Grant. 

\section{Triangulating the Cappell-Shaneson knot complement}\label{primarysec}

The triangulation of the figure-$8$ complement was discovered by a direct observation from a planar diagram of the figure-$8$ knot. The triangulation of the $2$-knot complement in this paper was discovered during an enumeration of $4$-dimensional triangulated manifolds.  As far as we are aware, there is as of yet no explicit description of Cappell-Shaneson knots embedded in $S^4$.  While in principle such descriptions exist \cite{Gom1, Gom2, Ak} in that many Cappell-Shaneson knots are known to be knots in the standard $S^4$, no explicit embeddings are known.  

Our terminology for triangulations is that a $0$-cell ($0$-simplex) is a {\it vertex}.  A $1$-cell ($1$-simplex) is an {\it edge}.  A $2$-cell ($2$-simplex) is a {\it triangle}.  A $3$-cell ($3$-simplex) is a {\it tetrahedron}.  A $4$-cell ($4$-simplex) is a {\it pentachoron}.  Plural of pentachoron will be {\it pentachora}. 

\begin{defn}\label{maintri} The central triangulation of this paper is described by gluing two disjoint ideal pentachora.  This means we take two disjoint $4$-simplices, and remove their vertices.  We then glue the ideal tetrahedra on the boundary according to the tetrahedron-pairing graph:
{
\psfrag{label}[tl][tl][1][0]{Figure 1}
\psfrag{a0}[tl][tl][0.8][0]{$0$} 
\psfrag{a1}[tl][tl][0.8][0]{$3$} 
\psfrag{a2}[tl][tl][0.8][0]{$2$} 
\psfrag{a3}[tl][tl][0.8][0]{$1$} 
\psfrag{a4}[tl][tl][0.8][0]{$4$} 
\psfrag{b0}[tl][tl][0.8][0]{$3$} 
\psfrag{b1}[tl][tl][0.8][0]{$0$} 
\psfrag{b2}[tl][tl][0.8][0]{$2$} 
\psfrag{b3}[tl][tl][0.8][0]{$1$} 
\psfrag{b4}[tl][tl][0.8][0]{$4$} 
\psfrag{g0}[tl][tl][0.8][0]{$0134 \leftarrow 1234$} 
\psfrag{g1}[tl][tl][0.8][0]{$0234 \leftarrow 0124$} 
\psfrag{g2}[tl][tl][0.8][0]{$0123 \leftarrow 0134$} 
\psfrag{h0}[tl][tl][0.8][0]{$0124 \to 1234$} 
\psfrag{g4}[tl][tl][0.8][0]{$0123 \to 0234$} 
$$\includegraphics[width=12cm]{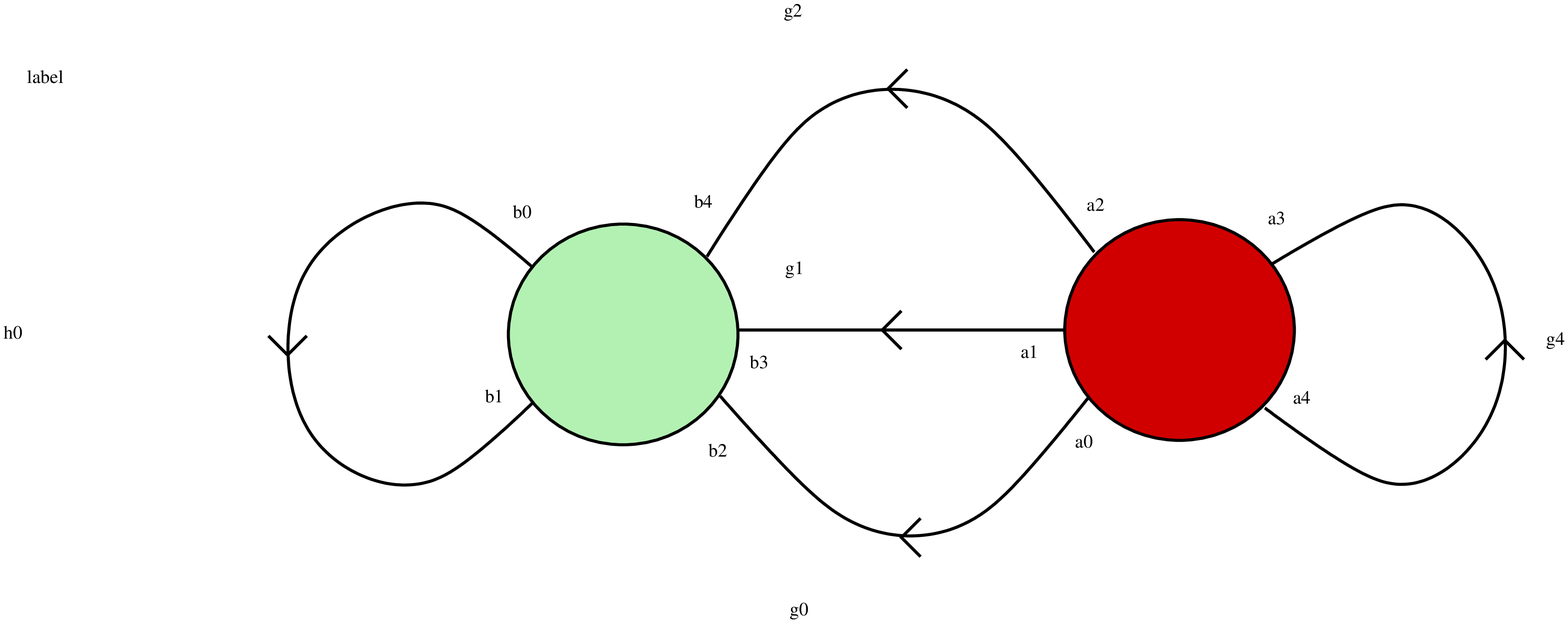}$$ 
}
What this picture means is the two pentachora are labelled `green' (left) and `red' (right) respectively.  The middle edge between the coloured circles emanates with label $3$ and terminates with $1$, this means the ideal tetrahedron opposite vertex $3$ in the red pentachoron is glued to the ideal tetrahedron opposite vertex $1$ in the green pentachoron.  The map of the ideal tetrahedron extends uniquely to an affine-linear map of the entire tetrahedron, which is prescribed by how it permutes the ideal vertices.  The ideal vertices of tetrahedron $3$ in the red pentachoron are sent via the map $(0,1,2,4) \longmapsto (0,2,3,4)$, that is, vertex $0$ is sent to $0$, $1$ sent to $2$, $2$ to $3$, and $4$ to $4$.  Similarly for the four other arrows.  This triangulation will be denoted $M$. 
\end{defn}

Notice, in the description of the triangulation $M$, we could have left-out the tetrahedron gluings such as $0124 \leftarrow 1234$, keeping only the labels such as $3$ and $0$ at the start and end of the arrows respectively -- since the gluings are order-preserving the gluing maps are prescribed by the index of the tetrahedra being glued together. 

\begin{lem}
$M$ is a triangulation of a manifold.
\begin{proof}
In general, if a triangulated $4$-manifold is specified by a tetrahedron-pairing graph as above, one must check the gluings do not induce self-identifications on the interior of any facet (edge, triangle, tetrahedron). Notice that all the gluing maps are order-preserving, so self-identifications are impossible.  Moreover, the triangulation $M$ admits an involution, which can be described as the unique map that switches the red and green pentachora, and reverses the labels of the vertices $(0,1,2,3,4) \longmapsto (4,3,2,1,0)$. 

If one runs through the gluings one quickly sees that all the (ideal) vertices are glued together, there is only one edge remaining after the gluings, $4$ triangles and $5$ tetrahedra.   To finish the proof, we argue that the link of the vertex is a triangulated $S^1 \times S^2$, and so this is called the `ideal boundary' or `cusp' and denoted $\partial M$.  This suffices to prove the triangulation is an ideal triangulation of a manifold.  The face-pairing graph for the link of the vertex is readily computed, in Figure 2.
{
\psfrag{a0}[tl][tl][0.8][0]{$3$} 
\psfrag{a1}[tl][tl][0.8][0]{$0$} 
\psfrag{a2}[tl][tl][0.8][0]{$2$}
\psfrag{a3}[tl][tl][0.8][0]{$1$}
\psfrag{a4}[tl][tl][0.8][0]{$4$}
\psfrag{g0}[tl][tl][0.8][0]{$3$} 
\psfrag{g1}[tl][tl][0.8][0]{$0$} 
\psfrag{g2}[tl][tl][0.8][0]{$1$}
\psfrag{g3}[tl][tl][0.8][0]{$2$}
\psfrag{g4}[tl][tl][0.8][0]{$4$}
\psfrag{0}[tl][tl][0.6][0]{$0$} 
\psfrag{1}[tl][tl][0.6][0]{$1$}
\psfrag{2}[tl][tl][0.6][0]{$2$}
\psfrag{3}[tl][tl][0.6][0]{$3$}
\psfrag{4}[tl][tl][0.6][0]{$4$}
$$\includegraphics[width=14cm]{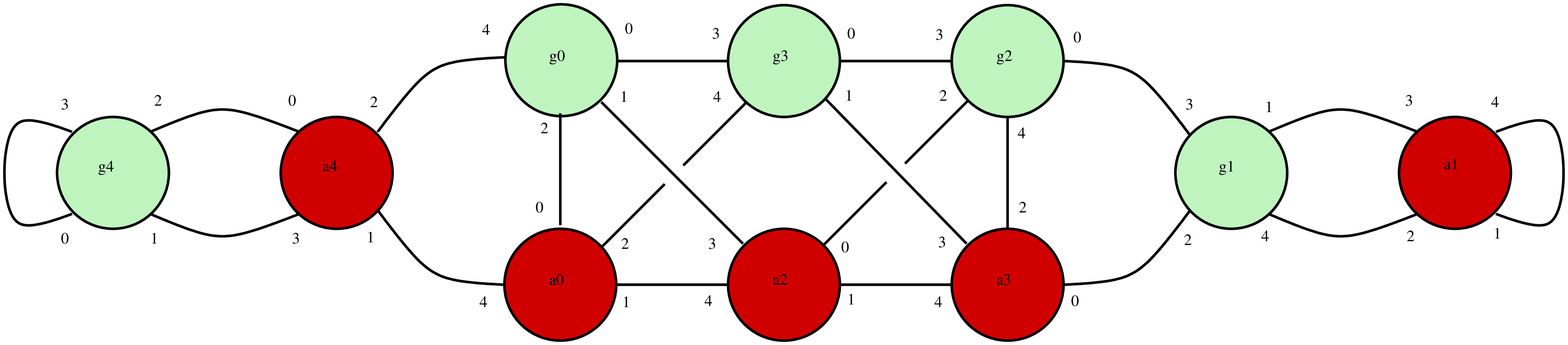}$$ 
\centerline{\it The face-pairing graph for the ideal boundary $\partial M$}
\centerline{Figure 2}
}

We will show that this triangulation is the union of two triangulated solid tori. In detail:  The leftmost tetrahedron has two of its faces glued together, giving a one-triangle triangulation of a M\"obius band.  The tetrahedron itself with those two faces glued together make a one-tetrahedron triangulation of $S^1 \times D^2$, called a {\it layered solid torus} \cite{JR}.  It's a quite beautiful and strange triangulation. The entire $1$-skeleton sits on $\partial (S^1 \times D^2)$, and only one $2$-cell is in the interior, being the triangulated M\"obius band. The union of the two tetrahedra (green and red, both labelled with $4$) comprise a $(2,3,-5)$-layered solid torus. This is a triangulation of $S^1 \times D^2$ such that there are exactly two triangles on the boundary.  Since there is a single $0$-cell on the boundary, every edge of the $1$-skeleton is a closed curve. Taking a boundary triangle one can coherently orient its boundary, and write down the signed intersection numbers with the meridian.  This is a well-defined triple (up to changing all the signs, and cyclic permutation). For this triangulation one gets $\pm (2,3,-5)$.  Attachment of either the green or the red tetrahedron labelled $3$ is along a single face -- which on the boundary amounts to a $1 \longmapsto 3$ Pachner move.  Attaching the other red / green tetrahedron labelled $3$ is a layering operation, in that on the boundary one is performing a $2 \longmapsto 2$ Pachner move.  Attaching either the green or red tetrahedron labelled $2$ is also a layering operation.   

Since the triangulation has the involution switching the red and green pentachoron, the ideal boundary also has such an involution, permuting the five green tetrahedral labelled $0,1,2,3,4$ with the five red labelled tetrahedra $4,3,2,1,0$ in that order.  So $\partial M$ is the union of two solid tori.  On the next page we will compute the dual polyhedral decomposition to the triangulation $M$, allowing us to compute $\pi_1 M$, and the homology of $M$.  From the homology long exact sequence of the pair $(\partial M, M)$ and Poincar\'e Duality, it follows that $H_1(\partial M)$ must be an infinite cyclic group, but $S^1 \times S^2$ is the only $3$-manifold with genus one Heegaard splitting and infinite $H_1$. 

It should be noted that the software Regina \cite{BBurton} can verify that a triangulated manifold is PL-equivalent to $S^1 \times S^2$, and this is how the authors first made this observation. 
\end{proof}
\end{lem}

We give a description of the dual polyhedral complex associated to the trianglation $M$.  It consists of:

\begin{itemize}
\item The two vertices at the barycentres of the pentachora, which we label $R$ and $G$ for red and green respectively.  
\item The five edges dual to the tetrahedra will be denoted $e_0, e_1, e_2, e_3, e_4$ and they correspond to $e_0 \equiv (0124\to 1234)$, $e_1 \equiv (0134\to 0123)$, $e_2 \equiv (0124\to 0234)$, $e_3 \equiv (1234\to 0134)$, $e_4 \equiv (0123\to 0234)$, where we give the orientation of the gluing map, from Definition \ref{maintri} respectively.  In particular, under our involution $e_i \longmapsto e_{4-i}^{-1}$ for all $i\in \{0,1,2,3,4\}$. 
\item The dual $2$-cells consist of two hexagons and two squares.  The hexagons will be denoted $f_1$ and $f_2$, the squares $f_3$ and $f_4$.  $f_1$ corresponds to the link about the $012$ face in the green pentachoron.  $f_2$ the link of $234$ in the red pentachoron. So the involution of $M$ permutes $f_1$ and $f_2$.  $f_3$ is the link of $234$ in the green pentachoron. $f_4$ the link of $012$ in the red pentachoron.  With appropriate choices of basepoints, a computation gives 
\begin{align*}
\partial f_1 = e_1e_4^{-1}e_2^{-1}e_3e_1^{-1}e_0 & \hskip 15mm \partial f_2 = e_3^{-1}e_0e_2e_1^{-1}e_3e_4^{-1} \\
\partial f_3 = e_2e_3^{-1}e_0^{-2} & \hskip 15mm \partial f_4 = e_2^{-1}e_1e_4^2.
\end{align*}
{
\psfrag{e0}[tl][tl][1][0]{$e_0$}
\psfrag{e1}[tl][tl][1][0]{$e_1$}
\psfrag{e2}[tl][tl][1][0]{$e_2$}
\psfrag{e3}[tl][tl][1][0]{$e_3$}
\psfrag{e4}[tl][tl][1][0]{$e_4$}
\psfrag{r0}[tl][tl][0.8][0]{$0$}
\psfrag{r1}[tl][tl][0.8][0]{$1$}
\psfrag{r2}[tl][tl][0.8][0]{$2$}
\psfrag{r3}[tl][tl][0.8][0]{$3$}
\psfrag{r4}[tl][tl][0.8][0]{$4$}
\psfrag{g0}[tl][tl][0.8][0]{$0$}
\psfrag{g1}[tl][tl][0.8][0]{$1$}
\psfrag{g2}[tl][tl][0.8][0]{$2$}
\psfrag{g3}[tl][tl][0.8][0]{$3$}
\psfrag{g4}[tl][tl][0.8][0]{$4$}
\centerline{Figure 3}
$$\includegraphics[width=14cm]{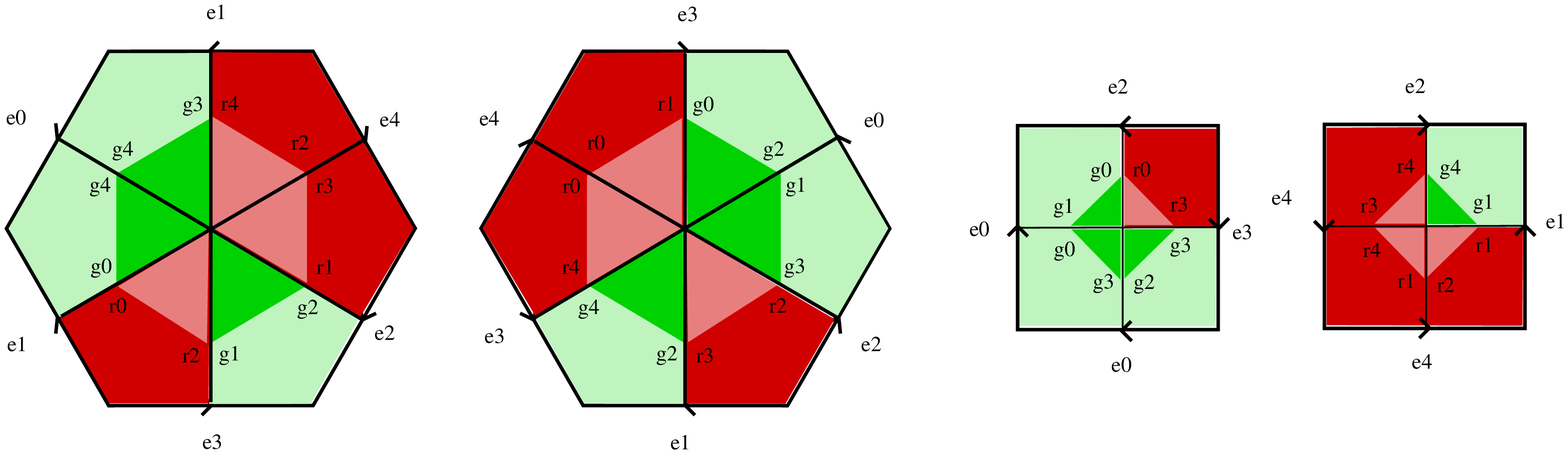}$$
\centerline{From left to right, the $2$-cells $f_1$, $f_2$, $f_3$ and $f_4$ respectively. }
}
\vskip 0.2cm
To make sense of the labellings, notice the periphery of these cells are labelled by the edges $\{e_0,e_1,e_2,e_3,e_4\}$, this indicates the attaching maps of the (dual) 2-skeleton to the (dual) $1$-skeleton.  We are traversing the boundary of the cells in Figure 3 in the counter-clockwise orientation, and our convention for concatenating paths is right-to-left.  The colours inside the cells indicate which pentachoron that part of the cell is in.  Each cell consists of a union of triangles, and these are the normal facets to the triangles of $M$ inside its pentachora.  So a red triangle with labels $2$ and $4$ (say, in $f_1$) indicates this is the edge normal to the triangle $013$ in the red pentachoron. 

\item There is only one $3$-cell. It is dual to the edge of $M$. 
\end{itemize}

The triangulation of the link of the edge in $M$ is included in Figure 4. Pictured are two combinatorial discs whose union is a sphere. Each disc includes some overlap with the other.  The discs are given as a union of red and green triangles.  At the three corners of the triangle are numbers.  For example, a green triangle with $4$, $1$ and $0$ in its corners means that it is the triangle normal to the edge $2,3$ in the green pentachoron.  Also in the picture is the dual polyhedral decomposition, showing how these normal triangles, when put together, can be expressed as a union of the dual $2$-cells $f_1, f_2, f_3$ and $f_4$.  This re-writing is more explicit in Figure 5. 
{
\psfrag{label}[tl][tl][1][0]{Figure 4}
\psfrag{r0}[tl][tl][0.6][0]{$3$}
\psfrag{r1}[tl][tl][0.6][0]{$0$}
\psfrag{r2}[tl][tl][0.6][0]{$2$}
\psfrag{r3}[tl][tl][0.6][0]{$1$}
\psfrag{r4}[tl][tl][0.6][0]{$4$}
\psfrag{g0}[tl][tl][0.6][0]{$3$}
\psfrag{g1}[tl][tl][0.6][0]{$0$}
\psfrag{g2}[tl][tl][0.6][0]{$1$}
\psfrag{g3}[tl][tl][0.6][0]{$2$}
\psfrag{g4}[tl][tl][0.6][0]{$4$}
$$\includegraphics[width=14cm]{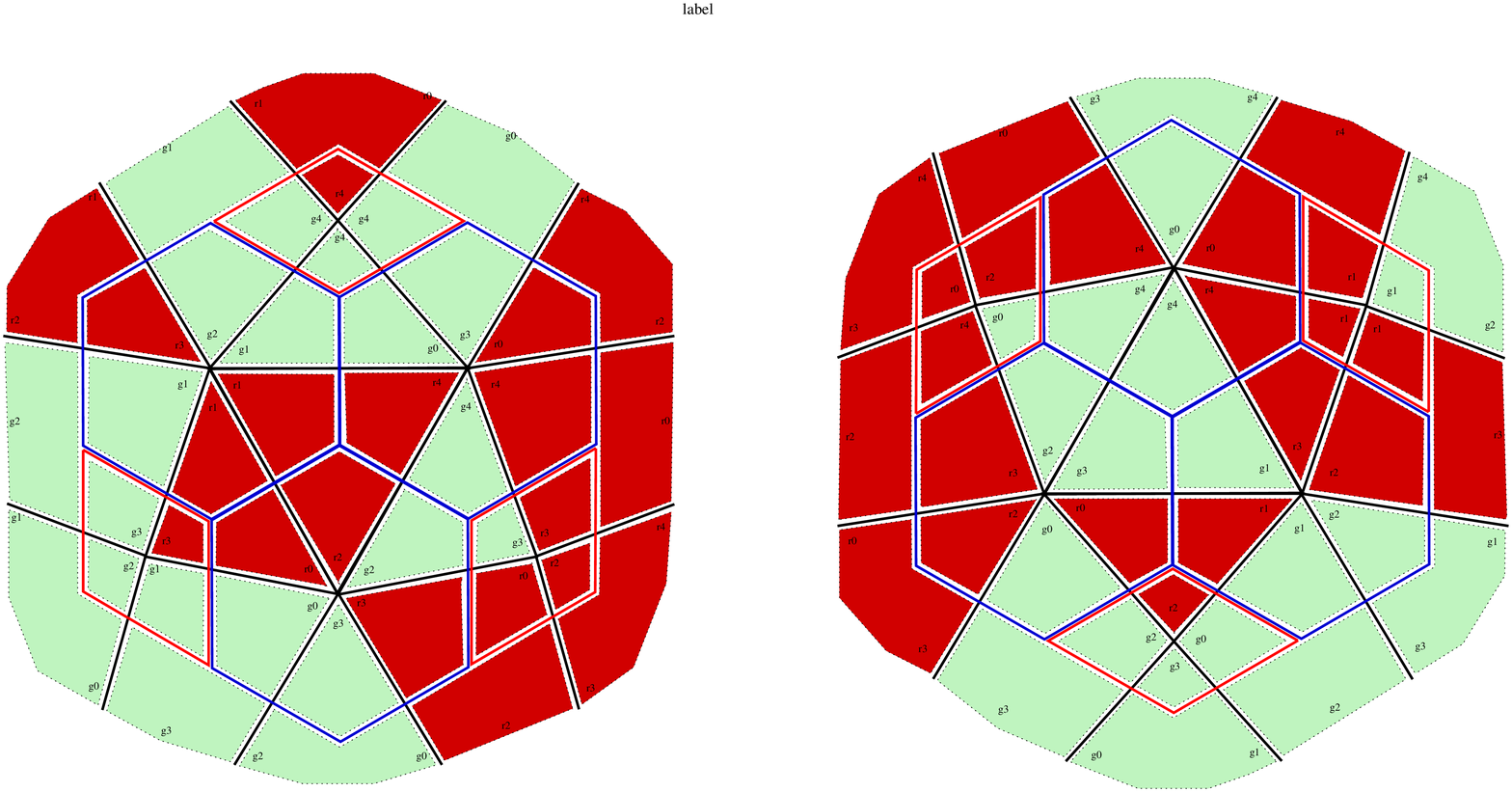}$$
\centerline{\it The link of the edge in $M$}
}

The attaching map of the triangulated edge link, given in the $2$-cells $\{f_1,f_2,f_3,f_4\}$ is in Figure 5. 
{
\psfrag{e1}[tl][tl][0.7][0]{$e_0$}
\psfrag{e2}[tl][tl][0.7][0]{$e_4$}
\psfrag{e3}[tl][tl][0.7][0]{$e_1$}
\psfrag{e4}[tl][tl][0.7][0]{$e_3$}
\psfrag{e5}[tl][tl][0.7][0]{$e_2$}
\psfrag{pf1}[tl][tl][0.8][0]{$+f_1$}
\psfrag{mf1}[tl][tl][0.8][0]{$-f_1$}
\psfrag{pf2}[tl][tl][0.8][0]{$+f_2$}
\psfrag{mf2}[tl][tl][0.8][0]{$-f_2$}
\psfrag{pf3}[tl][tl][0.8][0]{$+f_3$}
\psfrag{mf3}[tl][tl][0.8][0]{$-f_3$}
\psfrag{pf4}[tl][tl][0.8][0]{$+f_4$}
\psfrag{mf4}[tl][tl][0.8][0]{$-f_4$}

$$\includegraphics[width=10cm]{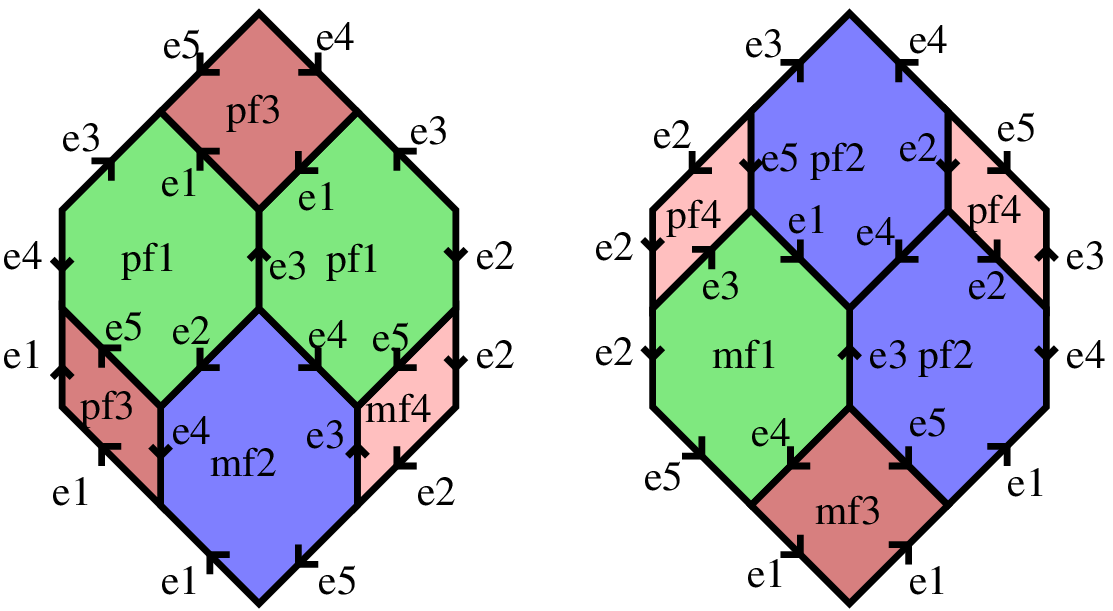}$$
\centerline{Figure 5}
}

\begin{prop}\label{csproc} The manifold $M$ triangulated in Definition \ref{maintri} is the complement of an embedded $S^2$ in a homotopy $S^4$. 
\begin{proof}
To verify this claim, we need to check that we can obtain a homotopy $4$-sphere by attaching a $D^2 \times S^2$ to the ideal boundary. 

As observed in the comments following Definition \ref{maintri}, the dual polyhedral decomposition has two $0$-cells $\{G,R\}$, five $1$-cells $\{e_0,e_1,e_2,e_3,e_4\}$, four two-cells $\{f_1,f_2,f_3,f_4\}$ and a three-cell.  If we collapse any one of $\{e_1,e_2,e_3\}$ the induced $1$-skeleton has a single $0$-cell, and so the induced $2$-skeleton gives us a presentation of the group. Collapsing $e_2$ gives us 
$$\langle e_0,e_1,e_3,e_4 | e_1e_4^{-1}e_3e_1^{-1}e_0, e_3^{-1}e_0e_1^{-1}e_3e_4^{-1}, e_3^{-1}e_0^{-2}, e_1e_4^2 \rangle$$
and we can use the last two relators to eliminate $e_1$ and $e_3$ respectively.  Re-labelling $e_0$ and $e_1$ to be $a$ and $b$ respectively gives us the presentation
$$ \pi_1 M \simeq \langle a, b | ba^2=a^3b^2, b^2a=a^2b^3 \rangle.$$
Notice that the  closed curve in $M$ corresponding to $a$ is parallel to a generator of $\pi_1 (\partial M)$, and if we were to append the relator $a=1$, the group would be trivial.  Since $M$ has the homology of a circle, attaching an $S^2 \times D^2$ to its boundary results in a homotopy $4$-sphere. 
\end{proof}
\end{prop}

The Cappell-Shaneson knots \cite{CS} are embeddings of $S^2$ in a family of homotopy $4$-spheres.  As of this date, it appears likely that all of these homotopy $4$-spheres are diffeomorphic to the standard smooth $S^4$ \cite{Ak, Gom1, Gom2}, although until quite recently these homotopy $4$-spheres were considered possible counter-examples to the smooth $4$-dimensional Poincar\'e Conjecture.   

Given a manifold $N$, and a diffeomorphism $f : N \to N$, there is an associated bundle $N \rtimes_f S^1$, this is `the bundle over $S^1$ with fibre $N$ and monodromy given by $f$.'  Its precise definition is the quotient space $N \times_{\Zed} \Real \equiv (N \times \Real)/\Zed$ where $\Zed$ acts on the product $N \times \Real$ by $n.(p,t) = (f^n(p),n+t)$, and where $f^n$ is the $n$-fold composite of $f$ with itself. Given a matrix $A \in GL(\Zed^3)$, by design $A$ acts on $\Real^3$, and also by design $A$ preserves $\Zed^3 \subset \Real^3$, thus $A$ acts on $\Real^3 / \Zed^3 \equiv (S^1)^3$.  Since $A$ fixes the origin in $\Real^3$, it fixes a point in $(S^1)^3$, call it $*$.  The Cappell-Shaneson manifold associated to $A$ is
$$ CS(A) = \left((S^1)^3 \setminus \{*\}\right) \rtimes_A S^1 $$
and $p_A(t)$ is defined to be the polynomial $p_A(t) = Det(A-tI)$. 

\begin{thm}\cite{CS} $CS(A)$ is the complement of a knot in a homotopy $4$-sphere provided $p_A(0)=1$, and $p_A(1) = \pm 1$. 
It is the complement of two inequivalent knots in homotopy $4$-spheres if further it is true that $p_A(t) > 0$ for all $t \in (-\infty,0) \subset \Real$. Moreover, $p_A(t)$ is the Alexander polynomial of the complement of these knots.
\end{thm}

It's also known that a Cappell-Shaneson manifold $CS(A)$ is the complement of two distinct knots in homotopy $4$-spheres if and only if $p_A(t) = 1+2t-t^2-t^3$, where `distinct' is in the sense of homeomorphisms of pairs $(S,K)$ where $S$ is a homotopy $4$-sphere and $K \subset S$ is a tame topologically embedded $2$-sphere.  See \cite{Hill} Theorem 18.5 for details.

\begin{prop}\label{pi1prop}
$\pi_1 M$ is isomorphic to the fundamental group of the complement of the Cappell-Shaneson manifold
$CS(A)$ with $A = \left[ \begin{array}{ccc} 0 & 0 & 1 \\
                        	 	    1 & 0 & 0 \\ 
                              		    0 & 1 & -1 \end{array} \right]$ with $p_A(t) = -t^3-t^2+1$.  
\begin{proof}
Think of the presentation
$$\pi_1 M \simeq \langle a, b | ba^2=a^3b^2, b^2a=a^2b^3 \rangle$$
as describing the $2$-skeleton of a CW-decomposition of $M$.  We lift the cell structure to the universal abelian cover $\widetilde{M}$ of $M$, to compute a presentation of $\pi_1 \widetilde{M}$.  This can also be thought of as the Reidemeister-Schreier technique, for the kernel of the abelianization map $\pi_1 M \to \Zed$.  Our abelianization map will send $a \longmapsto -1$, $b \longmapsto 1$.  Let $a_n$ and $b_n$ be the lifts of $a$ and $b$ respectively that emanate from the $n$-th lift of the $0$-cell, while $a_n$ terminates at the $(n-1)$-th lift, and $b_n$ at the $(n+1)$-th lift.   Crushing all the $b_n$'s to a point gives a wedge of circles.  This gives us the presentation
$$\pi_1 \widetilde{M} \simeq \langle \{ a_n : n \in \Zed\} | \{a_n^{-1}a_{n-1}^{-1}a_na_{n+1}a_{n+2}, a^{-1}_na_{n+2}a_{n+3} : n \in \Zed \} \rangle$$  
Using the relators $a_n^{-1}a_{n+2}a_{n+3}$ for all $n \in \Zed$, we can eliminate variables $a_n$ for $n \geq 3$ and $n < 0$, leaving only the generators $a_0, a_1, a_2$ and the relators $a_n^{-1}a_{n-1}^{-1}a_na_{n+1}a_{n+2}$, which once reduced via the relators $a^{-1}_na_{n+2}a_{n+3}$, gives the presentation
$$\pi_1 \widetilde{M} \simeq \langle a_0, a_1, a_2 | [a_0,a_1], [a_1,a_2], [a_0,a_2]\rangle.$$
The inclusion $\pi_1 \widetilde{M} \to \pi_1 M$ is given by
$a_i \longmapsto  b^{1-i}ab^i$.  The action of $\Zed$ on $\pi_1 \widetilde{M}$ can be thought of as being generated by conjugation by $b$, i.e. 
$z \longmapsto b^{-1}zb$, thus it is given by
$a_i \longmapsto a_{i+1}$ for $i=0,1$ and $a_2 \longmapsto a_2^{-1}a_0$.  The corresponding matrix $A$ is therefore
$$ A = \left[ \begin{array}{ccc} 0 & 0 & 1 \\
                                 1 & 0 & 0 \\ 
                                 0 & 1 & -1 \end{array} \right]. $$
\end{proof}
\end{prop}

\begin{prop}\cite{Hill} If a knot complement in a homotopy $4$-sphere has fundamental group isomorphic to the fundamental group of a Cappell-Shaneson manifold, then the knot complement is homeomorphic to the corresponding Cappell-Shaneson knot complement. 
\begin{proof}
This proposition is implicit in the work of Hillman \cite{Hill}.  We summarize Hillman's argument here. Given a knot complement $X$ in a homotopy $4$-sphere, we can obtain a closed $4$-manifold by attaching a $D^3 \times S^1$ to the boundary.   As an oriented smooth manifold, this is well-defined since all diffeomorphisms of $S^2 \times S^1$ extend to diffeomorphisms of $D^3 \times S^1$.  This manifold is called `the knot manifold' and we denote it $N$.  The inclusion $X \to N$ induces an isomorphism of fundamental groups, and $\chi(N) = \chi(X) + \chi (D^3 \times S^1) - \chi (S^2 \times S^1) = 0$.  To recover $X$ from $N$ one has to take the complement of an open tubular neighbourhood of the core $\{0\} \times S^1$ of the attached $D^3 \times S^1$ in $N$.  $\pi_1 X$ is normally generated by its meridional class, therefore the core $\{0\} \times S^1 \subset N$ normally generates $\pi_1 N$.  So the manifold $N$ together with this normal generator of $\pi_1 N$ determine $X$.  

Given a Cappell-Shaneson knot complement, the associated knot manifold is diffeomorphic to $(S^1)^3 \rtimes S^1$ where the monodromy is some orientation-preserving automorphism of $(S^1)^3$.  Moreover, up to conjugacy there are precisely two elements of the fundamental group which normally generate, one the reverse of the other.  Thus Cappell-Shaneson knot complements are determined by their knot manifolds.  If we start with a knot manifold $N$ whose fundamental group is isomorphic to that of a Cappell-Shaneson knot complement, Theorem 6.11 of \cite{Hill} states $N$ must be homeomorphic to a manifold of the form $(S^1)^3 \rtimes S^1$.  Therefore the corresponding complement is well-defined up to homeomorphism, so the manifold is homeomorphic to a Cappell-Shaneson knot complement.  Hillman's proof of Theorem 6.11 has two main steps.  Step one involves showing that $N$ is a $K(\pi,1)$, this is Corollary 3.5.1 of \cite{Hill}. Step two involves showing there is no obstruction to doing topological surgery, which follows from work of Farrell and Jones, plus Freedman's foundational work on topological surgery.  
\end{proof}
\end{prop}

\begin{cor}$M$ is homeomorphic to $CS(A)$.  
\end{cor}

\begin{prob}
Find a PL equivalence between $M$ and $CS(A)$. 
\end{prob}

Perhaps the most appealing way to approach this problem would be to find natural ideal triangulations of Cappell-Shaneson knot complements, perhaps in the spirit of Floyd and Hatcher's work \cite{FH} on bundles over $S^1$ whose fiber is a punctured torus. 

One of Thurston's broader observations in creating the ideal triangulation of the figure-$8$ knot complement is that knot and link complements in $S^3$ have fairly natural ideal polyhedral decompositions.  Given a link $L \subset S^3$, let $\nu L$ be an open tubular neighbourhood of $L$.  Further, let $\pi : S^3 \to [-1,1]$ be orthogonal projection onto some unit vector $v \in S^3$.  Isotope $L$ so that $\nu L \subset \pi^{-1}[-\frac{1}{2}, \frac{1}{2}]$.  Generically, $\pi : S^3 \setminus \nu L \to [-1,1]$ is a Morse function which has precisely one maximum, and one minimum, and whose only other critical points are on $\partial (\nu L)$ corresponding to local maximal and minima $\pi$ on $L$.  Stratified Morse Theory gives a cellular decomposition of the complement.  The corresponding presentation of $\pi_1 S^3 \setminus L$ is called a Wirtinger presentation of the link complement's fundamental group.  $\pi^{-1}(0)$ is a great $2$-sphere, so we can identify $\pi^{-1}[-\frac{1}{2},\frac{1}{2}]$ with $S^2 \times [-\frac{1}{2},\frac{1}{2}]$.  If orthogonal projection of $L$ onto $S^2$ is a regular link diagram, the link diagram induces a CW-decomposition of $S^2$ where the vertices are the crossings, and the edges are segments of the projected link.  This CW-structure lifts to a CW-decomposition of $S^3 \setminus L$ with only two $3$-cells corresponding to the max and min of $\pi$, together with a $2$-dimensional CW-complex which is the union of a CW-structure for $\partial (\nu L)$ together with $2$-dimensional cells that project to the CW-decomposition of $S^2$ corresponding to the planar link diagram.  There are also $1$-dimensional cells that project to the vertices of the link diagram.  

The upshot of this observation is that the complement of the link is the union of two $3$-balls along certain strata in their boundary.  An analogous construction builds polyhedral decompositions of link complements in $S^4$.  Given that this Cappell-Shaneson knot complement has such a simple triangulation, perhaps Cappell-Shaneson knots should appear quite early in any relatively small census of knots in $S^4$. 

\providecommand{\bysame}{\leavevmode\hbox to3em{\hrulefill}\thinspace}

\Addresses
\end{document}